\documentclass[12 pt]{article}
\usepackage{amsfonts}
\title{Competition between Discrete Random Variables, with Applications to Occupancy Problems}
\author{Julia Eaton\\
Department of Mathematics\\
University of Washington\\
Seattle, WA 98195, USA \and
Anant P.~Godbole\\
Department of Mathematics\\
East Tennessee State University\\
Johnson City, TN 37614, USA\and
Betsy Sinclair\\
Department of Political Science\\
University of Chicago\\
Chicago, IL 60637, USA}
\begin{document}
\def\ep{\varepsilon}
\def\lr{\left(}
\def\lc{\left\{}
\def\rc{\right\}}
\def\rr{\right)}
\def\nt{{n\choose2}}
\def\nr{{n\choose r}}
\def\ns{{n\choose s}}
\def\tv{{d_{{\rm TV}}}}
\def\P{{\rm Po}}
\def\cl{{\cal L}}
\def\ck{{\cal K}}
\def\p{\mathbb P}
\def\v{\mathbb V}
\def\e{\mathbb E}
\def\z{\mathbb Z}
\def\l{\lambda}
\newtheorem{thm}{Theorem}
\newtheorem{prop}[thm]{Proposition}
\maketitle
\begin{abstract}
Consider $n$ players whose ``scores" are independent and identically distributed values $\{X_i\}_{i=1}^n$ from some discrete distribution $F$.  We pay special attention to the cases where (i) $F$ is geometric with parameter $p\to0$ and (ii) $F$ is uniform on $\{1,2,\ldots,N\}$; the latter case clearly corresponds to the classical occupancy problem.  The quantities of interest to us are, first, the $U$-statistic $W$ which counts the number of ``ties" between pairs $i,j$; second, the univariate statistic $Y_r$, which counts the number of strict $r$-way ties between contestants, i.e., episodes of the form ${X_i}_1={X_i}_2=\ldots={X_i}_r$; $X_j\ne {X_i}_1;j\ne i_1,i_2,\ldots,i_r$; and, last but not least, the multivariate vector $Z_{AB}=(Y_A,Y_{A+1},\ldots,Y_B)$.  We provide Poisson approximations for the distributions of $W$, $Y_r$ and $Z_{AB}$ under some general conditions.  New results on the joint distribution of cell counts in the occupancy problem are derived as a corollary.
\end{abstract}
\section{Introduction}
\nocite{cmj}
\nocite{bhj}
\nocite{mori}
\nocite{bc}
In this paper we hope to shed new light on an old problem, studied extensively in, e.g., \cite{bhj}, \cite{kolchin}.  Consider $n$ players whose ``scores" are independent and identically distributed values $\{X_i\}_{i=1}^n$ from some discrete distribution $F$.  We consider the case of general distributions $F$ but pay special attention to the cases where (i) $F$ is geometric with parameter $p\to0$ and (ii) $F$ is uniform on $\{1,2,\ldots,N\}$; the latter case corresponds to the classical occupancy problem.  The quantities of interest to us are 
\begin{itemize}\item the $U$-statistic $W$ which counts the number of ``ties" between pairs $i,j$ (with $X_a=X_b=X_c=X_d$, for example, leading to a contribution of ${4\choose2}=6$ to the value of $W$); 
\item the univariate statistic $Y_r$ which counts the number of strict $r$-way ties between contestants, i.e., episodes of the form $X_i=x$ for some $x$ iff $i\in A$, $\vert A\vert=r$; and \item the multivariate vector $Z_{AB}=(Y_A,Y_{A+1},\ldots,Y_B)$.  
\end{itemize}
We provide Poisson approximations for the distributions of $W$, $Y_r$ and $Z_{AB}$ under some general conditions.  New results on the joint distribution of cell counts in the occupancy problem are derived as a corollary.

Consider the following elementary problem from \cite{cmj}: ``Two players use a coin that lands heads with probability $p$ to play a game that consists of a sequence of rounds.  In each round, the first player tosses the coin until a head appears.  Then the second player tosses the coin until a head appears.  If the players have the same number of flips in a round, the round is declared a tie and another round is played.  If not, the player with the larger number of flips wins the game.  Rounds are played successively until one of the two players wins the game."
Readers are asked  to find the expected number of rounds; the expected value of the total number of flips; and the probability distribution of the difference between the number of flips made by players 1 and 2 in a given round.  We briefly mention the solution for the first two of these questions:  The probability of a two person tie is clearly
\[ \sum\limits_{x=1}^\infty (1-p)^{2x-2} p^2 =
\frac{p}{2-p},\]
so that $\e(R)$, the expected number of rounds is given by 
\begin{eqnarray*}\e(R) &=&
\sum\limits_{x=1}^\infty x(\p({\rm tie}))^{x-1}(1-\p({\rm tie})) =
\sum\limits_{x=1}^\infty x\lr\frac{p}{2-p}\rr^{x-1}\lr1-\frac{p}{2-p}\rr\\ &=&
\frac{2-p}{2-2p},
\end{eqnarray*}
so that Wald's lemma yields for $\e(F)$, the expected total number of flips,
\[\e(F) = \e(
{F}/{R})\e(R)=\frac{2-p}{2-2p}\e(F/R)=\frac{2(2-p)}{p(2-2p)},\]
since the expected number $\e(F/R)$ of flips per round is clearly $2/p$.  
Computations for a three-person game, not mentioned in \cite{cmj}, are similar, but we need to lay down some rules as follows:  Three players each flip a $p$-coin until heads is flipped.  The player with the highest number of flips wins
unless there are ties between two or more players, in which case we repeat the process.  That is, the value of each of the three geometric variables in question must be unique.  We next compute the probability of a two- or three-way tie; the expected number of rounds; and the expected number of flips for $n=3$ -- to convince the reader that the situation rapidly becomes quite complicated as $n$ increases.  [The authors had a lively discussion with Lloyd Douglas, NSF Program Officer, about the following ``real-life" application of the $n$-person model with $p\to0$.  We wish to rank $n$ of the greatest free-throw shooters (or slam dunkers, or,...) in the National Basketball Association.  The players each shoot free throws until they miss -- conditional on the fact that no two players miss on the same attempt.  Rankings are then awarded in the obvious fashion.]

With three players (A,B,C), there are 3!=6 ways to have a strict
inequality and seven ways to tie, since
there is one way for a three way tie (which we loosely write as ``A = B = C") to
occur; ${3\choose 1}=3$
ways for A \(> \) B = C to occur; and
${3 \choose 2} = 3 $ ways for A =
B \(> \) C to occur.
Note that
 \begin{eqnarray*}\p(A = B = C) &=&p^3 + p^3(1-p)^3 + p^3(1-p)^6
 \ldots\\ &=& \sum\limits_{x=1}^\infty p^3(1-p)^{3x-3}\\ &=&
\frac{p^2}{3-3p+p^2}, \end{eqnarray*}
while the table below

\medskip

\begin{tabular}{lccc}
Case & A & B & C \\[.8ex]
 1 & TH, TTH, \(\ldots \) & H & H \\
 2 & TTH, TTTH, \(\ldots \) & TH & TH \\
 3 & TTTH, TTTTH, \(\ldots \) & TTH & TTH \\
 4 & TTTTH, TTTTTH, \(\ldots \) & TTTH & TTTH \\
 \vdots & \vdots & \vdots & \vdots
 \end{tabular}

\noindent reveals that
\[
\p(A>B=C)=
p^3 \sum\limits_{m=1}^\infty (1-p)^m \sum\limits_{i=0}^{m-1}
(1-p)^{2i} = \frac{p(1-p)}{3-3p+p^2}. 
\]
Finally, we observe from the table

\medskip

\begin{tabular}{lccc}
Case & A & B & C \\ [.8ex]
 1 & TH & TH & H \\
 2 & TTH & TTH & H, TH \\
 3 & TTTH & TTTH & H, TH, TTH \\
 4 & TTTTH & TTTTH & H, TH, TTH, TTTH \\
 \vdots & \vdots & \vdots & \vdots
 \end{tabular}

\noindent that 
\[\p(A = B>C) =p^3
\sum\limits_{m=1}^\infty(1-p)^{2m} \sum\limits_{i=0}^{m-1}(1-p)^i
= \frac{p(1-p)^2}{(2-p)(3-3p+p^2)},
\]
which leads to
\begin{eqnarray*}\p({\rm tie})&=&\p(A = B = C) + 3\p(A >B = C) + 3\p(A = B>C) \\
&=& \frac{5p^3 -
13p^2 + 9p}{(2-p)(3-3p+p^2)},
\end{eqnarray*}
and hence as before to
\[\e(R) = 
 \frac{1}{1-{\frac{5p^3 -13p^2 +
9p}{(2-p)(3-3p+p^2)}}}
\]
and
\[\e(F)=\e(F/R)\e(R)=\frac{3}{p}\lr\frac{1}{1-\frac{5p^3 - 13p^2 +
9p}{(2-p)(3-3p+p^2)}}\rr.\]
Competitions of the kind discussed above are best formulated in the more general context of occupancy models as follows:  $n$ balls are independently thrown into an infinite array of boxes so that any ball hits the $j$th box with probability $p_j$. Let $X_j$ be the number of balls in box $j$.  Then, with $p_j=(1-p)^{j-1}p$, we have the game inspired by \cite{cmj} ending iff $X_j\le1\ \forall j$.  Extremal versions of such questions have arisen in the literature before, often with surprising results.  Motivated by a question, posed by Carl Pomerance and arising in an additive number theory context, Athreya and Fidkowski \cite{af} proved that the probability $\pi_n$ that the highest numbered non-empty box has exactly one ball in it converges to a constant (which is shown to be one) iff $\lim_{n\to\infty}p_n/\sum_{j=n}^\infty p_j=0$. This is a condition that is not satisfied by, e.g., the sequence $p_n=1/2^n$ for which, quite interestingly, the limit superior and the limit inferior of the sequence $\pi_n$ differ in the {\it fourth} decimal place.  {These results had been independently obtained a few years earlier by Eisenberg et. al \cite{e2}, \cite {e1}, \cite{e3} and also by Bruss and O'Cinneide \cite{bc}.}  The comprehensive paper of M\'ori \cite{mori} is most relevant too: Here it is proven that given a double sequence of integer valued random variables, i.i.d.~within rows, and letting $\mu(n)$ denote the multiplicity of the maximal value in the $n$th row, the limiting distribution of $\mu(n)$ does not exist in the ordinary sense -- but that the intriguing empirical type a.s. limit result 
\[\lim_{t\to\infty}{{1}\over{\log t}}\sum_{n=1}^t{1\over n}I(\mu(n)=m)={{r^m}\over{m\log\lr{{1}\over{1-r}}\rr}},m=1,2,\ldots\]
holds, where $r$ is a parameter that depends on the distribution.  The whole field appears to be extraordinarily rich with known facts and tantalizing possibilities.

Results of the kind described above indicate that the cell counts for the $n$ person (Geometric) coin game are unlikely to behave in an asymptotically smooth way if $p=p_n\not\to0$.  This fact is borne out in Theorem 1, where we study the distribution of the number $W$ of pairs of equalities in the $n$ person game, with $W=0$ corresponding to the end of a ``round" in the sense of \cite {cmj}, and show that a good Poisson approximation is obtained if $np\to0$ (Geometric distribution) or $n/N\to0$ (Uniform distribution).  Theorem 2 concerns itself with the distribution of the number $Y_r$ of strict $r$-way ties (=the number of boxes with exactly  $r$ balls) and Theorem 3 with a multivariate generalization of Theorem 2.  The approximating distribution is Poisson (Theorem 2) or a product of independent Poisson variates (Theorem 3).  We note, moreover, that we were able to prove a result such as Theorem 3 relatively easily {\it probably} due to the approach taken -- we use as a counter the event that $r$ {\it specific} balls go into the same urn, rather than the conventional approach (e.g., \cite {bhj}, Section 6.2) of counting the number of urns with $r$ balls.  See also \cite{agg}, \cite{kolchin}, \cite{abt1},\cite{abt2}, and \cite{abt3}.
\section{Results}
\begin{thm}  Let $\{X_j\}_{j=1}^n$ be an integer valued sequence of i.i.d.~random variables with $\p(X_1=i)=p_i$, and consider the U-statistic
\[W=\sum_{K=1}^\nt I_K,\]
where, with $\ck$ denoting the $K$th 2-subset of $\{1,2,\ldots,n\}$, $I_K=1$ if $X_i=X_j; i,j\in\ck$ ($I_K=0$ otherwise).  For any two discrete random variables $T$ and $U$, let $\tv(\cl(T),\cl(U))$ denote the usual total variation distance between their distributions $\cl(T)$ and $\cl(U)$, i.e.,
\[\tv(\cl(T),\cl(U))=\sup_{A\subseteq\z^+}\vert\p(T\in A)-\p(U\in A)\vert,\]
and let $\P(\l)$ be the Poisson r.v.~with parameter $\l=\e(W)$.  Then
\[\tv(\cl(W),\P(\l))\le 2 n\pi+{{2 n\rho}\over{\pi}},\]
where $\pi=\p(X_1=X_2)$ and $\rho=\p(X_1=X_2=X_3)$.
\end{thm}

\medskip

\noindent{\bf Proof} The proof is an elementary application of, e.g., Theorem 2.C.5 in \cite{bhj}, which yields with $\l=\e(W)={n\choose2}\pi$
\begin{eqnarray*}&&\tv(\cl(W),\P(\l)\\
&&{}\le{{1-e^{-\l}}\over{\l}}\lr\sum_K\p^2(I_K=1)+\sum_K\sum_{\{L:\ck\cap\cl\ne\emptyset\}}\lc\e(I_KI_L)+\pi^2\rc\rr\\
&&{}\le\pi+{{2(n-2)\rho}\over{\pi}}+2(n-2)\pi\\
&&{}\le2n\pi+{{2n\rho}\over{\pi}},
\end{eqnarray*} 
as asserted.  

In Theorem 1, if the variables are uniform over $\{1,2,\ldots,N\}$, then $\pi=1/N; \rho=1/N^2$, so that $\tv(\cl(W),\P(\l))\le4n/N\to0$ if $N\gg n$, where, throughout this paper we write for $f_n,g_n\ge0$, $f_n\ll g_n$ (or $g_n\gg f_n$) if $f_n/g_n\to0$ as $n\to\infty$.  If the variables are Geometric($p$), then the discussion in Section 1 yields $\pi=p/(2-p)$ and $\rho=p^2/(3-3p+p^2)$, so that we get $\tv(\cl(W),\P(\l))\le 2np/(2-p)+2np(2-p)/(3-3p+p^2)\le6np\to0$ if $p\ll1/n$.  For the $n$ person game discussed in Section 1, we thus get
\begin{eqnarray*}
\p(W=0)&=&\p({\rm no\ ties})
=\exp\{-(n(n-1)p)/(2(2-p))\}\pm6np\\&=&\exp\{-\l\}\pm6np,
\end{eqnarray*}
\[\e(R)={{1}\over{\p(W=0)}}={{1}\over{e^{-\l}\pm6np}},\]
and
\[\e(F)={n\over p}\e(R).\]

The random variable $W$, while providing us with some insight, does not yield the level of detail that we desire.  For this reason, we turn our attention next to the variable $Y_r$ that counts the ``number of strict $r$-way ties," or, in other words, the ``number of boxes with exactly $r$ balls."  The development that follows is alternative to that provided, say, in \cite {bhj}, Theorems 6.C, 6.E, and particularly 6.F, though we do not make too many comparisons between our results and those of \cite{bhj}, since our main focus will be on the multivariate Theorem 3; the strategy of looking at specific sets of $r$ players is what sets our method apart.

Letting as before $\{X_j\}_{j=1}^n$ be an integer valued sequence of i.i.d.~random variables with $\p(X_1=i)=p_i$, we denote by
\[\pi=\sum_{x=1}^\infty p_x^r(1-p_x)^{n-r}\]
the probability that a specific set of $r$ players are involved in a strict tie, and thus
\[\l={n\choose r}\pi\]
is the expected number of boxes with exactly $r$ balls.  Throughout this paper we will employ, as in the previous sentence, the dual analogies of ``balls in boxes" and ``ties between contestants."  It may be readily verified that $\pi=(1/N^{r-1})(1-N^{-1})^{n-r}$ in the uniform case, and that in the geometric case $\pi=\sum_{x=1}^\infty(1-p)^{rx-r}p^r(1-(1-p)^{x-1}p)^{n-r}$ may be estimated as follows; the estimates may be seen to be tight provided that $p\to0$.  First we have
\begin{eqnarray*}
\pi&=&\sum_{x=1}^\infty(1-p)^{rx-r}p^r(1-(1-p)^{x-1}p)^{n-r}\\
&\ge&\sum\limits_{x=1}^\infty(1-p)^{rx-r}p^r[1-(n-r)(1-p)^{x-1}p] \\
&=&
\sum\limits_{x=1}^\infty(1-p)^{rx-r}p^r -
(n-r)\sum\limits_{x=1}^\infty(1-p)^{rx-r}p^r(1-p)^{x-1}p \\
&\ge&
\frac{p^{r-1}}{r} - \frac{2(n-r)p^r}{(r+1)(2-rp)},
\end{eqnarray*}
where the above inequalities follow since $(1-p)^r\ge 1-rp$ and $(1-p)^{r+1}\le 1-(r+1)p+[(r+1)rp^2]/2$.
Next note that
\begin{eqnarray*}
\pi&=&\sum_{x=1}^\infty(1-p)^{rx-r}p^r(1-(1-p)^{x-1}p)^{n-r}\\
&\le&\sum\limits_{x=1}^\infty(1-p)^{rx-r}p^r\cdot\\
{}&&\lr1-(n-r)(1-p)^{x-1}p +
\frac{(n-r)(n-r-1)}{2}(1-p)^{2x-2}p^2\rr\\
&=& \frac{p^r}{1-(1-p)^r} - \frac{(n-r)p^{r+1}}{1-(1-p)^{r+1}} +
\frac{(n-r)(n-r-1)}{2}\frac{p^{r+2}}{1-(1-p)^{r+2}}\\
&\le&
\frac{2p^{r-1}}{r(2-(r-1)p)} - \frac{(n-r)p^r}{r+1} +
\frac{(n-r)^2p^{r+1}}{(r+2)(2-(r+1)p)},
\end{eqnarray*}
so that in the geometric case, $\pi\sim p^{r-1}/r$ provided that $np^{(r+1)/2}\to0; rp\to0$.

We shall use the coupling approach as in \cite{bhj} to show that $\cl(Y_r)$ may be closely approximated by a Poisson distribution with the same mean.  We need to first find, given a sum $\sum_{j=1}^nI_j$ of indicator variables, a sequence $\{J_{ij}\}$ of indicator variables, defined on the same probability space as the $I_j$s, so that for each $j$,
\begin{equation}\cl(J_{1j},J_{2j},\ldots,J_{nj})=\cl(I_1,I_2,\ldots,I_n\vert I_j=1).\end{equation}
Good error bounds on a Poisson approximation are obtained if the $J_{ji}$s are chosen in a fashion that makes them ``not too far apart" from the $I_j$s.  We proceed in a manner similar to that in Theorem 6.F in \cite{bhj}, but the coupling we use is conditional, thus imparting a different flavor to the argument:  We have 
\[Y_r=\sum_{j=1}^\nr I_j,\]
where $I_j=1$ if the $j$th $r$-set is engaged in a strict tie.  Now we define the indicator variable $I_{jx}$ as being one if and only if $I_j$=1 and the members of the $j$th $r$-set all have ``value" $x$.  Now we proceed as follows:  If $I_{jx}=1$, we ``do nothing", setting $J_i=I_i$ for all $i$.  If, however, $I_{jx}=0$, we move all members of the $j$th $r$ set into the $x$th box (some of these might of course have occupied the $x$th box to begin with), while ejecting all its ``illegal" occupants and moving each these independently with probability $p_k/(1-p_x)$ to box $k;k\ne x$.  Finally we set $J_i=J_{ijx}=1$ if the $i$th $r$ set is involved in a strict tie {\it after} this interchange.  We need to verify that (1) holds in the modified form
\begin{equation}\cl(J_{1jx},J_{2jx},\ldots,J_{\nr jx})=\cl(I_1,I_2,\ldots,I_\nr\vert I_{jx}=1);\end{equation}
while this may be viewed as being ``obvious," we provide a proof next.  To show that (2) holds, it clearly suffices to show that any configuration (or sample point) corresponding to the members of the $j$th $r$-set being ``the only occupants of the $x$th box" is equally likely under both the conditional and unconditional models in (2). This strategy will achieve more, in fact, since we will not have to verify a condition similar to (2) when we move on to the multivariate case.  

 We let
 $a_i$ denote the score of the $i$th player not in
the $r$-clique in question ($a_i\ne x$), and 
 $b_i$  the score of the $i$th player in the $r$-clique, so that
$b_i = x$.  Now

\begin{eqnarray*}
\p({\rm configuration}\vert
 I_{jx}=1) 
  &=& {{\p({\rm configuration})}\over{\p(I_{jx}=1)}}\\
&=&\frac{\p(a_1, a_2, \ldots, a_{n-r},
b_{n-r+1}, \ldots, b_n)}{p_x^r(1-p_x)^{n-r}}
\\
&=&\frac{\p(a_1)\p(a_2)\ldots\p(a_{n-r})}{(1-p_x)^{n-r}}.
\end{eqnarray*}
Note also that the probability of the configuration under the coupled model is given by
\begin{equation}
\sum\limits_{l=0}^r {{r} \choose {l}}p_x^l[1-
p_x]^{r-l}\sum\limits_{S \subseteq {1, \ldots,n-r}} p_x^{\vert{S}\vert}\prod\limits_{j \in\{1,2,\ldots,n-r\}\setminus
S}p(a_j)\prod\limits_{j \in
S}\frac{p(a_j)}{1-p_x}. \end{equation}
Now 
\begin{eqnarray*}
{}&&\sum\limits_{S \subseteq {1, \ldots,n-r}}p_x^{\vert{S}\vert}\prod\limits_{j \in\{1,2,\ldots,n-r\}\setminus
S}p(a_j)\prod\limits_{j \in
S}\frac{p(a_j)}{1-p_x}\\ \quad&=&
\sum\limits_{t=0}^{n-r}{{n-r} \choose
{t}}\lr\frac{p_x}{1-p_x}\rr^t\prod\limits_{j=1}^{n-r}p(a_j)\\
&=& \lr\frac{p_x}{1-p_x} +
1\rr^{n-r}\prod\limits_{j=1}^{n-r}p(a_j) \\
 &=& \lr\frac{1}{1-p_x}\rr^{n-r}\prod_{j=1}^{n-r}p(a_j), \\
\end{eqnarray*}
which shows that (3) yields the same expression as before.  This proves the claim.

Now Theorem 2.B in \cite{bhj} leads to the following inequality:
\begin{eqnarray}
&&\tv(\cl(Y_r),\P(\l))\nonumber\\
&&{}\le\lr{{1-e^{-\l}}\over{\l}}\rr\sum_j\sum_x\p(I_{jx}=1)\lc\p(I_j=1)+\sum_{i\ne j}\p(I_i\ne J_{ijx})\rc\nonumber\\
&&{}\le\pi+\lr{{1-e^{-\l}}\over{\l}}\rr\sum_j\sum_x\p(I_{jx}=1)\sum_{i\ne j}\p(I_i\ne J_{ijx})
\end{eqnarray}
where $\l=\e(Y_r)$, $\pi=\p(I_j=1)$, and the coupled sequence $\{J_i\}=\{J_{ijx}\}$ satisfies (2) for each $j$ and $x$.
Consider first the case $\p(I_i=0,J_{ijx}=1)$, which is clearly impossible when $\vert i\cap j\vert\ge1$, and which we shall call Case I.  We thus have for $\vert i\cap j\vert=0$,
\begin{equation}
\sum_{i\ne j}\p(I_i=0,J_{ijx}=1)
={{n-r}\choose{r}}\sum_{y\ne x}\p(I_i=0,J_{ijx}=1,y),
\end{equation}
where the summand $\p(I_i=0,J_{ijx}=1,y)$ represents the probability that the $i$th $r$-set is not engaged in a strict $r$-way tie before the coupling, but is part of such a tie with common value $y$ after the coupling.  (5) thus yields 
\begin{eqnarray*}
\sum_{i\ne j}\p(I_i=0,J_{ijx}=1)&=&\sum_{i\ne j}\sum_{y\ne x}\lc\p(J_{ijx}=1,y)-\p(I_i=1,J_{ijx}=1,y)\rc\\&=&{{n-r}\choose{r}}\sum_{y\ne x}\mu_y
\end{eqnarray*}
where
\begin{eqnarray}
\mu_y&=&\sum_{q=0}^r{r\choose q}p_y^qp_x^{r-q}\sum_{s=0}^{n-2r}{{n-2r}\choose{s}}p_x^s(1-p_x-p_y)^{n-2r-s}\lr{{p_y}\over{1-p_x}}\rr^{r-q}\cdot\nonumber\\
&&{}\lr{{1-p_x-p_y}\over{1-p_x}}\rr^{s}\nonumber\\
&&{}-(1-p_y)^rp_y^r\sum_{s=0}^{n-2r}{{n-2r}\choose{s}}p_x^s(1-p_x-p_y)^{n-2r-s}\lr{{1-p_x-p_y}\over{1-p_x}}\rr^{s}\nonumber\\
&&{}=\sum_q{r\choose q}p_y^r\lr{{p_x}\over{1-p_x}}\rr^{r-q}\sum_s{{n-2r}\choose{s}}\lr{{p_x}\over{1-p_x}}\rr^{s}\cdot\nonumber\\&&{}(1-p_x-p_y)^{n-2r}\nonumber\\
&&{}-p_y^r(1-p_y)^r\sum_s{{n-2r}\choose{s}}\lr{{p_x}\over{1-p_x}}\rr^{s}(1-p_x-p_y)^{n-2r}\nonumber\\
&&{}=p_y^r\lr{{1}\over{1-p_x}}\rr^{n-2r}(1-p_x-p_y)^{n-2r}\lr{{1}\over{1-p_x}}\rr^{r}\nonumber\\
&&{}-p_y^r(1-p_y)^r(1-p_x-p_y)^{n-2r}\lr{{1}\over{1-p_x}}\rr^{n-2r}\nonumber\\
&&{}=p_y^r\lr1-{{p_y}\over{1-p_x}}\rr^{n-2r}\lc\lr{{1}\over{1-p_x}}\rr^{r}-(1-p_y)^r\rc
\end{eqnarray}
We now check to see the nature of the bound (6) in the uniform and geometric cases:  When the balls are distributed uniformly in $N$ boxes, we see that (6) leads to
\begin{eqnarray}
&&{}\sum_{i\ne j}\p(I_i=0,J_{ijx}=1)\nonumber\\&&{}={{n-r}\choose{r}}\sum_{y\ne x}\mu_y\nonumber\\
&&{}\le\nr\cdot N\cdot {{1}\over{N^r}}\lr1-{1\over N}\rr^{n-2r}\lc\lr{{N}\over{N-1}}\rr^r-\lr{{N-1}\over{N}}\rr^r\rc\nonumber\\
&&{}=\l \lr1-{1\over N}\rr^{-2r}\lc1-\lr1-{1\over N}\rr^{2r}\rc\nonumber\\
&&{}\le{{2r}\over{N}}\exp\{2r/(N-1)\}\l,
\end{eqnarray}
while in the geometric case we have 
\begin{eqnarray}
&&\sum_{i\ne j}\p(I_i=0,J_{ijx}=1)\nonumber\\
&&{}={{n-r}\choose{r}}\sum_{y\ne x}p_y^r\lr1-{{p_y}\over{1-p_x}}\rr^{n-2r}\lc{{1}\over{(1-p^x)^r}}-(1-p_y)^r\rc\nonumber\\
&&{}\le\nr\sum_{y=1}^\infty p_y^r(1-p_y)^{n-r}(1-p_y)^{-r}\lr(1-p_x)^{-r}-(1-p_y)^r\rr\nonumber\\
&&{}\le\l\max_{x,y}\lc(1-p_x)^{-r}(1-p_y)^{-r}-1\rc\nonumber\\
&&{}\le\l\lr{{1}\over{1-p}}\rr^{2r}\lc1-(1-p)^{2r}\rc\nonumber\\
&&{}\le\l\exp\{2rp/(1-p)\}2rp.
\end{eqnarray}
We now consider the case where $I_i=1$ and $J_{ijx}=0$ (Case II).  We clearly have $\p(I_i=1, J_{ijx}=0)=\p(I_i=1)$ for $\vert i\cap j\vert\ge1$ (Case II$'$), so we obtain
\begin{eqnarray}
\sum_{\vert i\cap j\vert\ge1}\p(I_i=1,J_{ijx}=0)&\le&r{{n-1}\choose{r-1}}\pi\nonumber\\
&=&{{r^2}\over{n}}\l.
\end{eqnarray} 
Next assume (Case II$''$) that $\vert i\cap j\vert=0$, and we seek to estimate the probability $\p(I_i=1, J_{ijx}=0)$.  If $y=x$, we bound $\p(I_i=1, J_{ijx}=0)$ by $\pi_x$, so that
\begin{equation}
\sum_{\vert i\cap j\vert=0}\p(I_{ix}=1, J_{ijx}=0)\le{{n-r}\choose{r}}\pi_x.
\end{equation}
If, on the other hand, $y\ne x$, then we have 
\begin{eqnarray*}
\sum_{\vert i\cap j\vert=0}\p(I_{i}=1, J_{ijx}=0)&=&{{n-r}\choose{r}}\sum_{y\ne x}\p(I_{iy}=1, J_{ijx}=0)\nonumber\\
&\le&{{n-r}\choose{r}}\sum_{y\ne x}\sum_{q\ge 1}{{n-2r}\choose{q}}p_x^qp_y^r{{qp_y}\over{1-p_x}};
\end{eqnarray*}
the above equation follows since in order for $I_{iy}=1, J_{ijx}=0$ to occur, we must have at least one of the $q$ ``bad" balls present in urn $x$ land in urn $y$ and thus ``spoil" the fact that $I_{iy}=1$.  We thus get
\begin{eqnarray}
\sum_{\vert i\cap j\vert=0}\p(I_{i}=1, J_{ijx}=0)&\le&{{n-r}\choose{r}}\sum_{y\ne x}{{p_y^{r+1}}\over{1-p_x}}\sum_{q\ge1}{{n-2r}\choose{q}}qp_x^q\nonumber\\
&\le&{{n-r}\choose{r}}n{{p_x}\over{1-p_x}}(1+p_x)^{n-2r-1}\sum_{y\ne x}{p_y^{r+1}}\nonumber\\
&\le&{{n-r}\choose{r}}n{{p_x}\over{1-p_x}}e^{np_x}\sum_{y\ne x}p_y^{r+1}.
\end{eqnarray}
Now (11) reduces in the uniform case to
\begin{eqnarray}
{{n-r}\choose{r}}{n\over {N-1}}e^{n/N}{{N-1}\over{N^{r+1}}}&\le&{{n}\over{N^2}}\l{{e^{n/N}}\over{\lr1-{1\over N}\rr^{n-r}}}\nonumber\\
&\le&{{n}\over{N^2}}\l{e^{n/N}}\exp\{(n-r)/(N-1)\}\nonumber\\
&\le&{{n}\over{N^2}}\l e^{2n/(N-1)}
\end{eqnarray}
and is bounded in the geometric case by
\begin{equation}
{{n-r}\choose{r}}{{np^{r+2}}\over{1-p}}\sum(1-p)^{(y-1)(r+1)}\le\l np^2e^{np}(1+o(1)),
\end{equation}
provided that $np^{(r+1)/2}\to0, rp\to0$.
Equations (4), (6), (9), (10) and (11) now yield
\begin{eqnarray}
&&\tv(\cl(Y_r),\P(\l))\nonumber\\
&&{}\le\pi+\lr{{1-e^{-\l}}\over{\l}}\rr\sum_j\sum_x\p(I_{jx}=1)\sum_{i\ne j}\p(I_i\ne J_{ijx})\nonumber\\
&&{}\le\pi+\lr1\wedge{1\over\l}\rr\sum_j\sum_x\p(I_{jx}=1)\bullet\nonumber\\
&&{}\bigg[{{n-r}\choose{r}}\sum_{y\ne x}p_y^r\lr1-{{p_y}\over{1-p_x}}\rr^{n-2r}\lc\lr{{1}\over{1-p_x}}\rr^{r}-(1-p_y)^r\rc+{{r^2}\over{n}}\l\nonumber\\
&&{}+{{n-r}\choose{r}}\pi_x+{{n-r}\choose{r}}n{{p_x}\over{1-p_x}}e^{np_x}\sum_{y\ne x}p_y^{r+1}\bigg].
\end{eqnarray}
We next evaluate (14) in the uniform case:  Equations (4), (7), (9), (10) and (12) give 
\begin{eqnarray}
&&\tv(\cl(Y_r),\P(\l))\nonumber\\
&&{}\le\pi+\lr1\wedge{1\over\l}\rr\sum_j\sum_x\p(I_{jx}=1)\sum_{i\ne j}\p(I_i\ne J_{ijx})\nonumber\\
&&{}\le \pi+\lr1\wedge{1\over\l}\rr\l\cdot\nonumber\\
&&{}\lc{{2r}\over{N}}\l\exp\{2r/(N-1)\}+{{r^2}\over{n}}\l+{{n-r}\choose{r}}{\pi\over N}+{{n}\over{N^2}}\l\exp\{2n/(N-1)\}\rc\nonumber\\
&&{}=\pi+\lr\l\wedge{\l^2}\rr\cdot\nonumber\\
&&{}\lc{{2r}\over{N}}\exp\{2r/(N-1)\}+{{r^2}\over{n}}+{1\over N}+{{n}\over{N^2}}\exp\{2n/(N-1)\}\rc.
\end{eqnarray}
We compare (15) with Equation 6.2.18 in \cite{bhj}, which yields the upper bound
\[\tv(\cl(Y_r),\P(\l))\le\lr\l\wedge{\l^2}\rr\lc{1\over N}+{{6n}\over{N^2}}+{{6r^2}\over{n}}\rc;\]
it is evident that (15) provides a better estimate if
\[{{2r}\over{N}}\exp\{2r/(N-1)\}\le{{5r^2}\over{n}}+(6-\exp\{2n/(N-1)\}){{n}\over{N^2}},\]
which is
a condition that holds under a wide range of circumstances, and certainly if $n/N\to0$.  Now in the geometric case, Equations (4), (8), (9), (10), and (13) reveal that (14) reduces as follows:
\begin{eqnarray}
&&\tv(\cl(Y_r),\P(\l))\nonumber\\
&&{}\le\pi+\lr1\wedge{1\over\l}\rr\sum_j\sum_x\p(I_{jx}=1)\bullet\nonumber\\
&&{}\lr 2\l rp\exp\{2rp/(1-p)\}+{{r^2}\over{n}}\l+\nr\pi_x+\l np^2\rr\nonumber\\
&&{}\le \pi+(\l\wedge\l^2)\lc2rp\exp\{2rp/(1-p)\}+{{r^2}\over{n}}+rp+np^2e^{np}\rc.
\end{eqnarray}
We have thus proved
\begin{thm}
Let $\{X_j\}_{j=1}^n$ be an integer valued sequence of i.i.d.~random variables with $\p(X_1=i)=p_i$.  Define $Y_r$ to be the number of strict $r$-way ties between these random variables.  Then the total variation distance between $\cl(Y_r)$ and a Poisson distribution with the same mean is given by (14).  This expression reduces to the one in Equation (15) when the distribution of the $X_i$s is uniform on $\{1,2,\ldots,N\}$ and to the expression in Equation (16) when $X_1\sim {\rm Geo}(p)$.
\end{thm}

\medskip

For the rest of the paper we will, for simplicity, restrict our attention to the classical occupancy problem of $n$ balls in $N$ boxes, assuming furthermore, that $n/N\to0$.  The goal is to obtain a multivariate Poisson approximation for the vector $Z_{AB}=\{Y_A,Y_{A+1}\ldots,Y_B\}$, for suitably restricted $A$ and $B$, and where the approximating Poisson vector consists of independent components.  First consider the quantities $\{\l_a:A\le a\le B\}.$  Since 
\begin{eqnarray*}\l_a&=&{n\choose a}\lr{1\over N}\rr^{a-1}\lr1-{1\over N}\rr^{n-a}\\&\sim&{{1}\over{\sqrt{2\pi a}}}N\lr{{ne}\over{Na}}\rr^a\exp\{(n-a)/N\}\\
&\sim&{{1}\over{\sqrt{2\pi a}}}N\lr{{ne}\over{Na}}\rr^a(1+o(1)),
\end{eqnarray*}
it follows, due to the fact that $n/N\to0$, that $\l_a$ is monotone decreasing in $a$.  Suppose that $\l_A<\infty$ for some finite $A$.  It then follows that the approximating Poisson distribution for $Y_{A+1}$ would have mean close to zero, making our agenda somewhat uninteresting.  We shall assume therefore that $\l_A\to\infty$ as $n,N\to\infty$.  Choices of the parameters that make this occur might be, e.g., $n=N^\alpha; \alpha<1$, when $\e(Y_a)\to0$ for all $a\ge A_0$, or, more interestingly, $n=N/\log N$ in which case the threshold $A_0$ would tend to infinity with $N$.  We thus seek values of $A$ and $B$ for which we get an ``interesting" multivariate Poisson approximation for the ensemble $(Y_A,\ldots,Y_B)$.  Now Theorem 10.J in \cite{bhj} yields, using notation suggested by that used in the proof of Theorem 2,
\begin{eqnarray}
&&{}\tv(\cl(Y_A,\ldots,Y_B),\prod_{a=A}^B\P(\l_a))\nonumber\\&&{}\le\sum_{a=A}^B\sum_{j=1}^{n\choose a}\sum_{x=1}^N\p(I_{ajx}=1)\lc\p(I_{aj}=1)+\sum_{biy\ne ajx}\p(I_{biy}\ne J_{biy})\rc,
\end{eqnarray}
where the last sum does not include the case $bi=aj$.  Correspondingly, we let $T_1,T_2,T_3$ denote the quantities
$$\sum_{a=A}^B\sum_{j=1}^{n\choose a}\sum_{x=1}^N\p(I_{ajx}=1)\p(I_{aj}=1)$$
$$\sum_{a=A}^B\sum_{j=1}^{n\choose a}\sum_{x=1}^N\p(I_{ajx}=1)\sum_{iy\ne jx}\p(I_{aiy}\ne J_{aiy})$$
and 
$$\sum_{a=A}^B\sum_{j=1}^{n\choose a}\sum_{x=1}^N\p(I_{ajx}=1)\sum_{b\ne a}\sum_{i=1}^{n\choose b}\sum_{y=1}^N\p(I_{biy}\ne J_{biy})$$ 
respectively; we need to compute the sum $T_1+T_2+T_3.$  First, we see that
\begin{eqnarray}
T_1&=&\sum_a\sum_j\p^2(I_{aj}=1)\nonumber\\
&=&\sum_a{n\choose a}\frac{1}{N^{2a-2}}\lr1-\frac{1}{N}\rr^{2n-2a}\nonumber\\
&\le&N^2\sum_{a=A}^B\lr\frac{ne}{aN^2}\rr^a\nonumber\\
&\le&N^2\sum_{a=2}^B\lr\frac{ne}{2N^2}\rr^a\nonumber\\
&=&\frac{e^2}{4}\frac{n^2}{N^2}(1+o(1))\to0
\end{eqnarray}
for each $A,B$.
The computation of $T_2$ follows as in the proof of Theorem 2.  The first component, $T_{21}$ is, by (7), given by
\begin{eqnarray}T_{21}&=&\sum_{a=A}^B\sum_j\sum_x\p(I_{ajx}=1)\cdot\frac{2a}{N}\exp\{2a/(N-1)\}\l_a\nonumber\\
&\le&\sum_a\l_a^2\frac{2a}{N}(1+o(1)).
\end{eqnarray}
Under what circumstances might the bound in (19) tend to zero?  Let us pause to consider this question before continuing.  If $n=\sqrt{N\log N}$ and $A=2$, then $\l_A\sim(\log N)/2$, the error bound of Theorem 2 is of magnitude $\sqrt{\log N/N}$, and the bound in (19) {\it does} approach zero.  However in this case $\l_3\to0$ so we are able to derive little useful beyond a Poisson approximation for $Y_2$, the number of ``days" with exactly two ``birthdays".  If $n=N^{0.9}$, then $\l_a\sim N^{1-0.1a}\to\infty$ for all $a=2,3,\ldots9$ but the summands in (19), asymptotically equal to $N^{1-0.2a}$, tend to zero only if $a=6,7,8,9$.  We thus have a potential multivariate approximation for $(Y_6,Y_7,Y_8,Y_9)$.  Finally, let $n=N/\log N$.  In this case, $\l_a\sim(e/a\log N)^a\cdot N$, and, with $a=\log N/(2\log\log N)$, for example, we see that
\begin{eqnarray*}
\l_a&\sim&\lr\frac{e}{a\log N}\rr^a\cdot N\\
&=&\lr\frac{2e\log \log N}{\log^2 N}\rr^{\frac{\log N}{2\log\log N}}\cdot \lr e^{2\log\log N}\rr^{\frac{\log N}{2\log\log N}}\\
&=&\lr{2e\log \log N}\rr^{\frac{\log N}{2\log\log N}}\to\infty,
\end{eqnarray*}
while with $a=\log N/[(4-\varepsilon)\log\log N]$ (we use $a=\log N/(3\log\log N)$ below) we have
\begin{eqnarray*}
\frac{\l_a^2}{N}&\sim&\lr\frac{e}{a\log N}\rr^{2a}\cdot N\\
&=&\lr\frac{3e\log \log N}{\log^2 N}\rr^{\frac{2\log N}{3\log\log N}}\cdot \lr e^{\frac{3\log\log N}{2}}\rr^{\frac{2\log N}{3\log\log N}}\\
&=&\lr{\frac{3e\log \log N}{\log^{1/2}N}}\rr^{\frac{2\log N}{3\log\log N}},
\end{eqnarray*}
which leads, with $A=\log N/(3\log\log N)$ and $B=\log N/(2\log\log N)$, to 
$$T_{21}\le 2(B-A)\frac{\l_A^2}{N}\le\frac{\log N}{3\log\log N}\lr{\frac{3e\log \log N}{\log^{1/2}N}}\rr^{\frac{2\log N}{3\log\log N}}\to0;$$
we thus have a potential Poisson approximation for the vector $(Y_A\ldots Y_B)$.  Next note that the term $T_{22}$ that corresponds to (9) is given by
\begin{equation}
T_{22}=\sum_{a=A}^B\sum_{j=1}^{n\choose a}\sum_{x=1}^N\p(I_{ajx}=1)\frac{a^2}{n}\l_a=\sum_a\frac{a^2}{n}\l_a^2.
\end{equation}
Finally we combine the two remaining terms (10) and (12), as reflected in (15), to get
\begin{eqnarray}
T_{23}&=&\sum_{a=A}^B\sum_{j=1}^{n\choose a}\sum_{x=1}^N\p(I_{ajx}=1)\lr\frac{{n\choose a}\pi_a}{N}+\frac{n}{N^2}\l_a\exp\{2n/(N-1)\}\rr\nonumber\\
&=&\sum_a\frac{\l_a^2}{N}+\frac{n}{N^2}\l_a^2(1+o(1)).
\end{eqnarray}
Turning to a computation of $T_3$, we first observe that for $a\ne b$ and $\vert i\cap j\vert\ge1$, it is impossible for $I_{biy}=0,J_{biy}=1$ to occur.  Accordingly, as in the calculation leading up to (6) we see that
\begin{eqnarray}
&&{}\sum_b\sum_i\sum_y\p(I_{biy}=0,J_{biy}=1)\nonumber\\
&&{}=\sum_b\sum_i\sum_y\{\p(J_{biy}=1)-\p(I_{biy}=1,J_{biy}=1)\}\nonumber\\
&&{}=\sum_b{{n-a}\choose{b}}\sum_y\Huge\{\sum_{q=0}^b{b\choose q}\lr\frac{1}{N}\rr^q\lr\frac{1}{N}\rr^{b-q}\cdot\nonumber\\
&&{}\sum_s{{n-a-b}\choose{s}}\lr\frac{1}{N}\rr^s\lr1-\frac{2}{N}\rr^{n-a-b-s}\lr\frac{1}{N-1}\rr^{b-q}\lr\frac{N-2}{N-1}\rr^s-\nonumber\\
&&{}\lr\frac{N-1}{N}\rr^a \lr\frac{1}{N}\rr^b\sum_s{{n-a-b}\choose{s}}\lr\frac{1}{N}\rr^s\lr1-\frac{2}{N}\rr^{n-a-b-s}\lr\frac{N-2}{N-1}\rr^s\Huge\}\nonumber\\
&&{}=\sum_b{{n-a}\choose{b}}N\cdot\nonumber\\
&&{}\lc\lr\frac{1}{N-1}\rr^b\lr\frac{N-2}{N-1}\rr^{n-a-b}-\lr\frac{N-1}{N}\rr^a\lr\frac{1}{N}\rr^b\lr\frac{N-2}{N-1}\rr^{n-a-b}\rc\nonumber\\
&&{}\le\sum_b\l_b\frac{(a+b)}{N}(1+o(1)).
\end{eqnarray}
As in the univariate case, $\p(I_{biy}=1,J_{biy}=0)=\p(I_{biy}=1)$ if $\vert i\cap j\vert\ge1$.  Hence
\begin{eqnarray}
\sum_b\sum_{\vert i\cap j\vert\ge1}\sum_y\p(I_{biy}=1,J_{biy}=0)&\le&\sum_ba{{n-1}\choose{b-1}}\pi_b\nonumber\\
&=&\sum_b \frac{ab}{n}\l_b.
\end{eqnarray}
If, however, $\vert i\cap j\vert=0$, then
\begin{equation}
\sum_b\sum_{\vert i\cap j\vert=0}\p(I_{bix}=1,J_{bix}=0)\le\sum_b{{n-a}\choose{b}}\pi_{bx},
\end{equation}
and, being rather crude with the final estimation
\begin{eqnarray}
&&{}\sum_b\sum_{\vert i\cap j\vert=0}\sum_{y\ne x}\p(I_{biy}=1,J_{biy}=0)\nonumber\\
&&{}\le\sum_b{{n-a}\choose{b}}\frac{1}{N^b}\sum_y\sum_{q\ge1}{{n-a-b}\choose{q}}\frac{1}{N^q}\frac{q}{N-1}\nonumber\\
&&{}\le\sum_b{n\choose b}\frac{n}{N^{b+1}}\sum_{q\ge1}{{n-1}\choose{q-1}}\frac{1}{N^{q-1}}(1+o(1))\nonumber\\
&&{}\le\sum_b\frac{n}{N^2}\l_b(1+o(1)).
\end{eqnarray}
Collecting equations (18) through (25), we see that the following holds:
\begin{thm}  When $n$ balls are randomly assigned to $N$ boxes, where $n\ll N$, the joint distribution of the multivariate vector $(Y_A,\ldots,Y_B)$ of exact box counts may be approximated by a Poisson vector with independent components.  More specifically,
$$\tv\lr\cl(Y_A,\ldots,Y_B),\prod_{a=A}^B\P(\l_a) \rr\le \varepsilon_{n,N,A,B}$$
where $\l_a=\e(Y_a)={n\choose a}\frac{1}{N^{a-1}}\lr1-\frac{1}{N}\rr^{n-a}$ and $\varepsilon_{n,N,A,B}$ is of magnitude
$$\sum_{a=A}^B\l_a^2\lr\frac{2a}{N}+\frac{a^2}{n}+\frac{1}{N}+\frac{n}{N^2}\rr+\sum_a\l_a\lr\sum_{b\ne a}\l_b\lr \frac{(a+b)}{N}+\frac{ab}{n}+\frac{1}{N}+\frac{n}{N^2}\rr\rr.$$  In addition, an application, e.g., of Theorem 10.K in \cite{bhj} may provide slight improvements in the above, through a partial reinstatement of the so-called ``magic factor".
\end{thm}

\medskip

\noindent{\bf Acknowledgment}  The research of all three authors was supported by NSF Grants DMS-0049015 and DMS-0552730, and was conducted at East Tennessee State University, when Eaton and Sinclair were undergraduate students at Rochester University and the University of the Redlands respectively.  An earlier version of this work was presented at IWAP-Piraeus and a completed version was presented at Lattice Path Combinatorics-Johnson City.

\end{document}